\documentclass[11pt,french]{amsart}
\usepackage[T1]{fontenc}
\usepackage[latin1]{inputenc}
\usepackage{latexsym}
\usepackage{amsthm}
\usepackage{amssymb}

 \theoremstyle{plain}
\newtheorem{thm}{Th\'eor\`eme}[section]
  \theoremstyle{d\'efinition}
  \newtheorem{defn}[thm]{D\'efinition}
  \theoremstyle{plain}
  \newtheorem{prop}[thm]{Proposition}
  \theoremstyle{d\'efinition}
  \newtheorem{condition}[thm]{Hypoth\`ese}
  \theoremstyle{plain}
  \newtheorem{lem}[thm]{Lemme}
  \theoremstyle{plain}
  \newtheorem{cor}[thm]{Corollaire}
  \theoremstyle{remarque}
  \newtheorem{notation}[thm]{Notation}
  \theoremstyle{remarque}
  \newtheorem{rem}[thm]{Remarque}
 \theoremstyle{d\'efinition}
  \newtheorem{example}[thm]{Exemple}

\AtBeginDocument{
  
}

\usepackage{babel}
\addto\extrasfrench{\providecommand{\fg}{\ifdim\lastskip>\z@\unskip\fi~\frqq}}

\begin{document}

\title{L'application cotangente des surfaces de type général}

\author{Xavier Roulleau}
\begin{abstract}
We study surfaces of general type $S$ whose cotangent sheaf is generated
by its global sections. We define a map called the cotangent map of
$S$ that enables us to understand the obstructions to the ampleness
of the cotangent sheaf of $S$. These obstructions are curves on $S$
that we call {}``non-ample''. We classify surfaces with an infinite
number of non-ample curves and we partly classify the non-ample curves. 
\end{abstract}
\maketitle


\section*{Introduction.}

On étudie dans le présent article l'amplitude du fibré cotangent des
surfaces de type général quand ce fibré cotangent est engendré par
sections globales.\\
 Les propriétés des variétés à fibré cotangent ample sont discutées
en \cite{Debarre2}. A l'aide d'une construction due à Bogomolov,
Conduché et Palmieri \cite{Conduche} ont montré que les ratios $\frac{c_{1}^{2}}{c_{2}}$
des nombres de Chern de telles surfaces forment un sous-ensemble dense
dans l'intervalle $[1,2]$. Une lecture attentive de leur papier et
l'utilisation du théorème des sections hyperplanes de Lefschetz permettent
de raffiner leur résultat comme suit:\\
Pour tout entier $q\geq12$, la fermeture des ratio de Chern des
surfaces d'irrégularité $q$, à fibré cotangent ample et engendré
par ses sections globales, contient un intervalle non-vide $I_{q}\subset[1,2]$.
La réunion des intervalles croissants $I_{q}$ est l'intervalle $]1,2[$.

Plutôt que de construire des surfaces dont le fibré cotangent est
ample, on propose d'étudier les obstruction à son amplitude. Ces obstructions
bien comprise, on peut alors espérer caractériser les surfaces à fibré
cotangent ample.

Rappelons d'abord la définition d'amplitude d'un fibré vectoriel $\mathcal{E}$
sur une variété $X$, telle qu'introduite par Hartshorne :\\
Notons $\mathbb{P}(\mathcal{E}^{*})$ le projectivisé du dual de
$\mathcal{E}$, $\pi:\mathbb{P}(\mathcal{E}^{*})\rightarrow X$ la
projection et $\mathcal{O}_{\mathbb{P}(\mathcal{E}^{*})}(1)$ le fibré
inversible tautologique tel que :\[
\pi_{*}\mathcal{O}_{\mathbb{P}(\mathcal{E}^{*})}(1)=\mathcal{E}.\]

\begin{defn}
Le fibré $\mathcal{E}$ est dit ample si $\mathcal{O}_{\mathbb{P}(\mathcal{E}^{*})}(1)$
est ample.
\end{defn}
Quand le fibré vectoriel $\mathcal{E}$ est engendré par l'espace
de ses sections globales, on dispose du critère d'amplitude suivant,
dû à Gieseker \cite{Gieseker}:
\begin{prop}
\label{proposition:critere d'amplitude}(Gieseker) Le fibré $\mathcal{E}$
est ample si et seulement si pour toute courbe $C\hookrightarrow X$,
le fibré $\mathcal{E}\otimes\mathcal{O}_{C}$ n'a pas de quotient
isomorphe à $\mathcal{O}_{C}$. 
\end{prop}
Pour aborder la question de l'amplitude du fibré cotangent, nous sommes
amenés à poser l'hypothèse suivante sur la surface considérée:
\begin{condition}
\label{conditions}Nous travaillerons dans ce qui suit avec une surface
$S$ définie sur $\mathbb{C}$, lisse, de type général, dont le fibré
cotangent $\Omega_{S}$ est engendré par l'espace de ses sections
globales $H^{o}(\Omega_{S})$ et d'irrégularité $q=\dim H^{o}(\Omega_{S})$
vérifiant $q>3$.
\end{condition}
L'article \cite{Conduche} de Conduché et Palmieri montre que de telles
surfaces existent abondamment. Notons $T_{S}=\Omega_{S}^{*}$ le fibré
tangent, $\pi:\mathbb{P}(T_{S})\rightarrow S$ la projection et $\mathcal{O}_{\mathbb{P}(T_{S})}(1)$
le fibré tautologique. On dispose d'une identification naturelle :\[
H^{o}(\mathbb{P}(T_{S}),\mathcal{O}_{\mathbb{P}(T_{S})}(1))\simeq H^{o}(\Omega_{S}).\]

\begin{defn}
\label{definition appli cot}Par hypothèse, le morphisme naturel $H^{o}(\Omega_{S})\otimes\mathcal{O}_{S}\rightarrow\Omega_{S}$
est surjectif, le morphisme $H^{o}(\Omega_{S})\otimes\mathcal{O}_{\mathbb{P}(T_{S})}\rightarrow\mathcal{O}_{\mathbb{P}(T_{S})}(1)$
est donc également surjectif et définit un morphisme : \[
\psi:\mathbb{P}(T_{S})\rightarrow\mathbb{P}(H^{o}(\Omega_{S})^{*})=\mathbb{P}^{q-1}\]
appelé l'\textbf{application cotangente} de la surface. 
\end{defn}
Ce morphisme et la question de l'amplidude du fibré cotangent sont
les principaux objets d'étude du présent article. En vue de la proposition
\ref{proposition:critere d'amplitude}, on pose la définition suivante
\begin{defn}
\label{definition: courbe non-ample}Une courbe $C\hookrightarrow S$
est dite \textbf{non-ample} si et seulement si le fibré $\Omega_{S}\otimes\mathcal{O}_{C}$
possède un quotient isomorphe à $\mathcal{O}_{C}$. 
\end{defn}
Les courbes non-amples constituent donc l'obstruction à l'amplitude
du fibré cotangent ; nous en étudions les propriétés. L'application
cotangente permet de traduire les propriétés algébriques du fibré
cotangent par les propriétés géométriques de l'image de $\psi$. Ainsi
une courbe $C\hookrightarrow S$ est non-ample si et seulement si
il existe une section $t:C\hookrightarrow\mathbb{P}(T_{S})$ contractée
en un point $p$ par $\psi$. En ce cas, l'image de $\pi^{*}C$ par
$\psi$ est un cône de sommet $p$.

Il est donc légitime de poser la définition suivante:
\begin{defn}
\label{definition:point  exceptionnel}Un point $p$ de $\mathbb{P}^{q-1}$
est dit exceptionnel si la fibre de $\psi$ en $p$ est de dimension
$>0$.
\end{defn}
Nous noterons $\Delta$ le lieu des points exceptionnels. Un problème
naturel est d'étudier $\Delta$ et de caractériser les surfaces pour
lesquelles ce fermé est de dimension strictement positive.

Résumons les résultats obtenus en vue de ces problèmes:
\begin{prop}
L'image de l'application cotangente est de dimension $3$. \\
Soit $p$ un point de $\mathbb{P}^{q-1}$. La fibre $\psi^{-1}(p)$
est de dimension au plus $1$. \\
Le lieu des points exceptionnels est de dimension au plus $1$. 
\end{prop}
Nous établissons ensuite une borne sur le degré de l'application cotangente
et sur le degré de son image en fonction des nombres de Chern $c_{1}^{2},c_{2}$
de la surface. 

Soit $C\hookrightarrow S$ une courbe, on appellera suite cotangente
de $C$ la suite exacte : \[
0\rightarrow\mathcal{O}_{C}(-C)\rightarrow\Omega_{S}\otimes\mathcal{O}_{C}\rightarrow\Omega_{C}\rightarrow0.\]
 Le théorème suivant classifie partiellement les courbes non-amples
de la surface :
\begin{thm}
\label{kodaira2}Soit $C$ une courbe de la surface $S$, alors :\\
a) La courbe $C$ est non-ample et vérifie $C^{2}<0$ si et seulement
si $C$ est lisse de genre $1$.\\
b) La courbe $C$ est non-ample et vérifie $C^{2}=0$ si et seulement
si $C$ est lisse de genre $>1$ et la suite cotangente est scindée. 
\end{thm}
Ce théorème s'obtient à partir d'un critère de lissité de Lipman \cite{lipman}.
La propriété a) nous semble un résultat particulièrement intéressant
: étant donné une courbe $C$ vérifiant $C^{2}<0$ sur une surface
lisse $X$, il est facile de construire une fibré vectoriel $\mathcal{E}$
sur $X$ engendré par ses sections globales et tel que $C\hookrightarrow S$
soit la seule courbe pour laquelle $\mathcal{E}\otimes\mathcal{O}_{C}$
admette un quotient trivial. Ceci illustre de nouveau les proprietés
particulières que possède le fibré cotangent d'une surface.

Rappelons qu'une surface fibrée $f:S\rightarrow B$ est dite isotriviale
si ses fibres lisses sont isomorphes entre-elles. Donnons des exemples
de courbes non-amples sur des surfaces:
\begin{prop}
a) Les fibres lisses d'une fibration de $S$ sont non-amples si et
seulement si la fibration est isotriviale.\\
b) Soit $S=C^{(2)}$ le produit symétrique d'une courbe $C$ de
genre $>3$ et non-hyperelliptique. Cette surface vérifie les hypothèses
\ref{conditions} et contient une infinité de courbes non-amples $C$
vérifiant $C^{2}=1$. Son image est la variété des sécantes d'une
courbe.
\end{prop}
La partie a) résulte d'un critère dû à Martin-Deschamps \cite{Martin}.\\
On dispose ainsi d'exemples de surfaces ayant une infinité de courbes
non-amples. Nous classifions ces surfaces dans le théorème suivant: 
\begin{thm}
\label{classification des surfaces ayant une infinit=0000E9 de courbes non-amples}Soit
$S$ une surface possédant une infinité de courbes non-amples et telle
que le lieu exceptionnel $\Delta$ ne soit pas une droite de $\mathbb{P}^{q-1}$.
L'image de l'application cotangente vérifie alors l'une des deux propriétés
suivantes :\\
a) La surface $S$ est une surface fibrée isotriviale dont les
fibres sont les courbes non-amples. Le lieu des points exceptionnels
$\Delta$ est formé de deux courbes et l'image de l'application cotangente
est la variété développée par les sécantes de ces deux courbes.\\
b) L'image de l'application cotangente est la variété des sécantes
d'une courbe. Une courbe non-ample $C\hookrightarrow S$ vérifie en
ce cas : $C^{2}>0$. Il existe un morphisme fini de $S$ dans le produit
symetrique $B^{(2)}$ d'une courbe $B$.
\end{thm}
Je remercie vivement Igor Reider qui m'a proposé ce thème de travail.
Je tiens a remercier également le Max-Planck Institute de Bonn où
une partie de cet article a été rédigé.

\section{Etude de l'application cotangente.}

\subsection{Les définitions de l'application cotangente.\label{sub:Les-d=0000E9finitions-de}}

Rappelons quelques proprietés partagées par une surface $S$, $A$
sa variété d'Albanese, $\vartheta:S\rightarrow A$ un morphisme d'Albanese
et $\Omega_{S}$ son fibré cotangent :
\begin{lem}
\label{lemme voisin}(\cite{griffithsharris} p. 331). La différentielle
$d\vartheta:T_{S}\rightarrow\vartheta^{*}T_{A}=H^{o}(\Omega_{S})^{*}\otimes\mathcal{O}_{S}$
du morphisme $\vartheta$ est le dual du morphisme d'évaluation $H^{o}(\Omega_{S})\otimes\mathcal{O}_{S}\rightarrow\Omega_{S}$.
\end{lem}
Notons $p_{r}:S\times H^{o}(\Omega_{S})^{*}\rightarrow H^{o}(\Omega_{S})^{*}$
la projection sur le second facteur. Par définition, le morphisme
$\vartheta$ est une immersion locale si le morphisme : \[
p_{r}\circ d\vartheta_{s}:T_{S,s}\rightarrow H^{o}(\Omega_{S})^{*}\]
est injectif en tout point $s$ de $S$. Le lemme \ref{lemme voisin}
permet ainsi une interprétation géométrique de la condition \ref{conditions}:
\begin{cor}
\label{alba}Le morphisme d'Albanese d'une surface est une immersion
locale si et seulement si le fibré cotangent de la surface est engendré
par ses sections globales.\end{cor}
\begin{notation}
Maintenant et pour le reste de cet article, $S$ est une surface vérifiant
les hypothèses \ref{conditions} et on conservera les notations de
l'introduction.
\end{notation}
Le corollaire suivant est une conséquence immédiate du lemme \ref{lemme voisin}
et permet une interprétation plus géométrique de l'application cotangente:
\begin{cor}
\label{retrouver s}L'application cotangente de $S$ est le projectivisé
du morphisme : \[
p_{r}\circ d\vartheta:T_{S}\rightarrow H^{o}(\Omega_{S})^{*}.\]

\end{cor}
L'application cotangente est donc le morphisme qui à un point de la
surface et à une direction tangente associe la direction tangente
dans la variété d'Albanese.

Regardons maintenant les propriétés de la restriction de l'application
cotangente en la fibre $\pi^{-1}(s)$ d'un point $s$. La fibre en
$s$ de la projection $\pi$ est la courbe $\mathbb{P}(T_{S,s})\simeq\mathbb{P}^{1}$
et la restriction de $\mathcal{O}_{\mathbb{P}(T_{S})}(1)$ à cette
courbe est le fibré de degré $1$, donc: 
\begin{lem}
\label{pi est injectif sur fibre de cot}La restriction de l'application
cotangente à la fibre $\pi^{-1}(s)$ ($s$ point de $S$) est un plongement
et son image est une droite de $\mathbb{P}^{q-1}$.\end{lem}
\begin{notation}
\label{notation de la droite Ls}Pour un point $s$ de $S$, on notera
$L_{s}\hookrightarrow\mathbb{P}^{q-1}$ la droite image de la fibre
$\pi^{-1}(s)$ par $\psi$.
\end{notation}
L'image de l'application cotangente est donc la réunion des droites
$L_{s}$ ($s$ point de $S$). Notons $G(2,q)$ la grassmannienne
des sous-espaces vectoriels de dimension $2$ de $H^{o}(\Omega_{S})^{*}$.
Cette grassmannienne paramètre également les droites de l'espace projectif
$\mathbb{P}^{q-1}$.
\begin{defn}
\label{d=0000E9finition du morphisme de gauss}Le \textbf{morphisme
de Gauss} : \[
\mathcal{G}:S\rightarrow G(2,q)\]
de la surface $S$ est défini par la surjection $H^{o}(\Omega_{S})\otimes\mathcal{O}_{S}\rightarrow\Omega_{S}$. 
\end{defn}
Par construction, le point $\mathcal{G}(s)$ représente la droite
projective $L_{s}$, ou encore :
\begin{cor}
\label{diff=0000E9rentielle de l'apli. cot.} Le morphisme de Gauss
est le morphisme qui à un point $s$ de la surface associe le point
de $G(2,H^{o}(\Omega_{S})^{*})$ représentant le plan : \[
p_{r}\circ d\vartheta_{s}(T_{S,s})\subset H^{o}(\Omega_{S})^{*}.\]

\end{cor}
Le corollaire $2$ de \cite{Ran} montre que :
\begin{lem}
\label{gauss est fini}(Ran) Le morphisme de Gauss est fini sur son
image.
\end{lem}
Soit $U$ le fibré universel de $G(2,q)$ et $\mathbb{P}(U)$ le projectivisé
de $U$. La situation étudiée est donc la suivante : \[
\begin{array}{ccccc}
\mathbb{P}(T_{S}) & \stackrel{\widetilde{\psi}}{\rightarrow} & \mathbb{P}(U) & \stackrel{\pi_{1}}{\rightarrow} & \mathbb{P}^{q-1}\\
\pi\downarrow &  & \downarrow\pi_{2}\\
S & \stackrel{\mathcal{G}}{\rightarrow} & G(2,q)\end{array}\]
où $\pi_{1},\pi_{2}$ sont les projections naturelles et le morphisme
$\widetilde{\psi}$ vérifie $\pi_{1}\circ\widetilde{\psi}=\psi$.

\subsection{Dimension des fibres et de l'image de l'application cotangente.}

Cette section porte sur les fibres de l'application cotangente et
la dimension de son image. Nous établissons que ces fibres sont de
dimension au plus $1$ et qu'il y a au plus une famille de dimension
$1$ de fibres $\psi^{-1}(p)$ de dimension $1$. Nous montrons ensuite
que l'image de l'application cotangente est de dimension $3$. 
\begin{lem}
\label{pi est injectif sur la fibre de appli cot}Soit $p$ un point
de $\mathbb{P}^{q-1}$. Le morphisme $\pi$ est injectif sur les points
de la fibre $\psi^{-1}(p)$.\end{lem}
\begin{proof}
Soit $s$ un point de $S$. Si l'intersection de $\pi^{-1}(s)$ et
de $\psi^{-1}(p)$ est non vide, cette intersection est nécessairement
un point car $\psi$ est un plongement sur la fibre $\pi^{-1}(s)$
(cf. lemme \ref{pi est injectif sur fibre de cot}).
\end{proof}
Notons $F$ l'image de l'application cotangente. Si $p$ est un point
de $\mathbb{P}^{q-1}$, nous noterons $D_{p}$ l'image de la fibre
$\psi^{-1}(p)$ par $\pi$ : le fermé sous-jacent est formé des points
$s$ de la surface tels que la droite $L_{s}$ passe par $p$. 
\begin{prop}
\label{pi est injectif}Pour tout point $p$ de $\mathbb{P}^{q-1}$,
la fibre $\psi^{-1}(p)$ et le schéma $D_{p}$ sont de dimension au
plus $1$.\end{prop}
\begin{proof}
L'application cotangente $\psi$ n'est pas constante donc pour tout
point $p$ de $\mathbb{P}^{q-1}$, la dimension de $\psi^{-1}(p)$
est inférieure ou égale à $2$.\\
Supposons que la fibre de $\psi$ en un point $p$ soit de dimension
$2$. Soit $S'$ le schéma réduit associé à la fibre $\psi^{-1}(p)$.
Par le lemme \ref{pi est injectif sur la fibre de appli cot}, la
restriction $\pi_{|S'}$ du morphisme $\pi$ à $S'$ est injective
sur les points de $S'$. Puisque $\pi$ est propre, le morphisme :
\[
\pi_{|S'}:S'\rightarrow S\]
 est surjectif. Le morphisme $\pi_{|S'}$ est bijectif, séparable
dans la variété normale $S$ : c'est un isomorphisme (cf. \cite{milne}
remarque 6.21). Il existe donc un morphisme $t:S\rightarrow\mathbb{P}(T_{S})$
tel que $\pi\circ t:S\rightarrow S$ soit l'identité. A ce morphisme
$t$ correspond un quotient : \[
\Omega_{S}\rightarrow\mathcal{L}\]
 où $\mathcal{L}$ est le fibré inversible qui vérifie $t^{*}\mathcal{O}_{\mathbb{P}(T_{S})}(1)=\mathcal{L}$
( \cite{hart} chapitre II proposition 7.12). Le fibré $\mathcal{O}_{\mathbb{P}^{q-1}}(1)$
est trivial au point $p$ et, puisque $\psi\circ t$ contracte $S$
au point $p$, le fibré : \[
t^{*}\psi^{*}\mathcal{O}_{\mathbb{P}^{q-1}}(1)=\mathcal{L}\]
est trivial. Ainsi $\Omega_{S}$ possède un quotient d'image un fibré
inversible trivial. Puisque $\Omega_{S}$ est engendré par $H^{0}(\Omega_{S})$,
ce quotient a une section et $\Omega_{S}=\mathcal{O}_{S}\oplus\mathcal{O}_{S}(K)$
($K$ diviseur canonique). Cela est impossible car nous avons supposé
$S$ de type général. La fibre de $\psi$ en un point est donc de
dimension au plus $1$.\\
Pour tout point $p$ de $\mathbb{P}^{q-1}$, l'image de la fibre
$\psi^{-1}(p)$ par $\pi$ est donc également de dimension inférieure
ou égale à $1$.
\end{proof}
Soit $p$ un point de $\mathbb{P}^{q-1}$ tel que $D_{p}=\pi(\psi^{-1}(p))$
soit de dimension $1$. Considérons $D$ une composante irréductible
de dimension $1$ de $D_{p}$ munie de sa structure réduite, alors
:
\begin{prop}
\label{p sommet d'un cone} L'image par l'application cotangente de
la surface réglée $\pi^{-1}(D)$ est un cône de sommet $p$. \end{prop}
\begin{proof}
L'application $\psi$ est injective sur les fibres de $\pi$, le fermé
irréductible $\psi(\pi^{-1}(D))$ est donc au moins de dimension $1$.
\\
Si $\psi(\pi^{-1}(D))$ est de dimension $1$, alors c'est une
droite projective $L\hookrightarrow\mathbb{P}^{q-1}$ et pour tout
point $s$ de $D$, on a : $L_{s}=L$. Mais cela est impossible car
le morphisme de Gauss $\mathcal{G}$ est fini (lemme \ref{gauss est fini}).\\
Donc $\psi(\pi^{-1}(D))$ est une surface et puisque toutes les
droites $L_{s}$ ($s$ point de $D$) passent par le point $p$, c'est
un cône de sommet $p$.\end{proof}
\begin{defn}
Un point $p$ de $\mathbb{P}^{q-1}$ est dit exceptionnel si la fibre
$\psi^{-1}(p)$ est de dimension $1$. \label{notation points exceptionnels}
\end{defn}
Soit $\Delta$ l'ensemble des points exceptionnels. Par la proposition
\ref{p sommet d'un cone}, chaque point de $\Delta$ est le sommet
d'un cône. Etudions la géométrie de $\Delta$ : 
\begin{prop}
\label{f est de dim 3}L'image de l'application cotangente est de
dimension $3$ ; $\Delta$ est un fermé de $\mathbb{P}^{q-1}$ vide
ou de dimension inférieure ou égale à $1$.\end{prop}
\begin{proof}
Notons $F$ l'image de l'application cotangente et $\psi_{|F}$ la
restriction de $\psi$ à son image. \\
Si les fibres générique de $\psi_{|F}$ sont de dimension $1$
alors $F$ est de dimension $2$. En ce cas, pour tout point $p$
de $F$, le schéma $D_{p}$ est de dimension $1$ et le fermé $\psi(\pi^{-1}(D_{p,red}))$
est de dimension $2$ contenu dans la surface irréductible $F$. Cela
implique que $\psi(\pi^{-1}(D_{p,red}))$ est égal à $F$. Ainsi deux
points quelconques de $F$ sont sommets de cônes et sont reliés par
une droite $L$ contenue dans $F$ pour laquelle il existe $s$ tel
que $L=L_{s}$.\\
On en déduit que $F$ est un plan projectif. Or la variété $F$
est non-dégénérée dans $\mathbb{P}^{q-1}$, ainsi : $q=3$. Mais on
a fait l'hypothèse que la surface $S$ est d'irrégularité $q>3$.
\\
Ainsi, pour un point $p$ générique de l'image de $\psi$, le schéma
$D_{p}$ est de dimension nulle et $\psi$ est génériquement finie
sur son image.\\
Le morphisme $\psi$ étant propre et génériquement fini, l'ensemble
des points exceptionnels $\Delta$ est un sous-schéma fermé de $\mathbb{P}^{q-1}$
(cf. \cite{hart} p.94).
\end{proof}

\subsection{Degré de l'application $\psi$, degré de $\psi_{*}\pi^{*}C$.}

\subsubsection{Degré de l'application cotangente et de son image.\label{paragraphe degr=0000E9 appli cot}\protect \\
}

Après avoir étudié les fibres de l'application cotangente, nous étudions
son degré.

Soit $F\hookrightarrow\mathbb{P}^{q-1}$ l'image de l'application
cotangente, notons $\deg F$ son degré et $\deg\psi$ le degré de
l'application $\psi$ sur son image. 
\begin{prop}
\label{degre de psi}Soit $c_{1}$ et $c_{2}$ les classes de Chern
de la surface, alors : \[
\deg F\deg\psi=c_{1}^{2}[S]-c_{2}[S].\]
\end{prop}
\begin{proof}
Soit $H$ une section hyperplane de $F$, notons : $h=\psi^{*}H$.
Par définition des classes de Chern de $\Omega_{S}$, le cycle $h$
vérifie : \[
h^{2}-h.\pi^{*}c_{1}(\Omega_{S})+\pi^{*}c_{2}(\Omega_{S})=0.\]
Puisque $c_{1}=-c_{1}(\Omega_{S})$ et $c_{2}=c_{2}(\Omega_{S})$,
on a donc : $h^{3}=-h^{2}.\pi^{*}c_{1}-h.\pi^{*}c_{2}$ et par substitution,
on obtient :\[
h^{3}=h.\pi^{*}c_{1}^{2}+\pi^{*}c_{2}.\pi^{*}c_{1}-h.\pi^{*}c_{2}=h.\pi^{*}(c_{1}^{2}-c_{2}).\]
Le degré du terme de droite est égal à celui de $c_{1}^{2}-c_{2}$.
De plus : $\deg h^{3}=\deg H^{3}\deg\psi$ et $\deg H^{3}=\deg F$,
donc $\deg F\deg\psi=c_{1}^{2}[S]-c_{2}[S]$ où $c_{1}^{2}[S]$ et
$c_{2}[S]$ sont les nombres de Chern de la surface.
\end{proof}
Rappelons le diagramme suivant, noté ({*}): \[
\begin{array}{ccccc}
\mathbb{P}(T_{S}) & \stackrel{\widetilde{\psi}}{\rightarrow} & \mathbb{P}(U) & \stackrel{\pi_{1}}{\rightarrow} & \mathbb{P}^{q-1}\\
\pi\downarrow &  & \downarrow\pi_{2}\\
S & \stackrel{\mathcal{G}}{\rightarrow} & G(2,q)\end{array}\]
Le morphisme $\psi$ vérifie : $\pi_{1}\circ\widetilde{\psi}=\psi$.
Notons $\pi_{1}'$ la restriction de $\pi_{1}$ à l'image de $\widetilde{\psi}$
sur l'image de $\widetilde{\psi}\circ\pi_{1}$ ainsi : 
\begin{prop}
\label{degr=0000E9 de psi divise c1^2}Le degré de $\psi$ est le
produit du degré du morphisme de Gauss $\mathcal{G}$ par le degré
de $\pi'_{1}$.\end{prop}
\begin{rem}
L'image de l'application cotangente étant non-dégénérée, son degré
est supérieur ou égal à $q-3$, donc le degré de la restriction de
l'application cotangente $\psi_{|}:\mathbb{P}(T_{S})\rightarrow F$
est majoré par $\frac{c_{1}^{2}[S]-c_{2}[S]}{q-3}$.
\end{rem}

\subsubsection{Degré du morphisme de Gauss.\protect \\
}

Comme nous allons le voir, on peut créer des surfaces $X$ dont le
degré du morphisme de Gauss $\mathcal{G}_{X}$ est arbitraire. Nous
saisissons l'occasion pour construire également des surfaces dont
le fibré cotangent n'est pas ample:\\
La donnée d'un revêtement est un triplet $\Lambda=(D,\mathcal{L},m)$
où $D$ est un diviseur lisse sur $S$, $\mathcal{L}$ un fibré inversible
et $m>1$ un entier tel que $\mathcal{L}^{m}=\mathcal{O}_{S}(D)$.
Si $D=0$, on demande de plus que $m$ soit le plus petit entier $k>0$
tel que $\mathcal{L}^{k}=\mathcal{O}_{S}$. A une telle donnée est
associé un revêtement cyclique $\tau_{\Lambda}:X\rightarrow S$ ramifié
en $D'=\tau_{\Lambda}^{-1}D$ où $X=X(\Lambda)$ est une surface lisse.
On a :
\begin{prop}
\label{pro:a)-revetement etales}A) Si la surface $S$ ne possède
pas de fibrations, alors il existe une infinité de données de revêtements
$\Lambda=(0,\mathcal{L},m)$ tels que:\\
i) la surface $X=X(\Lambda)$ vérifie l'hypothèse \ref{conditions},
\\
ii) le morphisme $\tau_{\Lambda}^{*}:H^{0}(\Omega_{S})\rightarrow H^{0}(\Omega_{X})$
est un isomorphisme qui permet d'identifier ces deux espaces. \\
iii) sous cette identification, on a : $\mathcal{G}_{X}=\mathcal{G}_{S}\circ\tau_{\Lambda}$.
En particulier le degré de $\mathcal{G}_{X}$ est divisible par $m$
et l'application cotangente de $X$ a même image que celle de $S$.
Ainsi $\Omega_{X}$ est ample si et seulement si $\Omega_{S}$ est
ample.\\
B) Si $D$ est ample, alors le fibré cotangent de $X$ n'est pas
engendré par ses sections globales au-dessus de $D'$ et n'est pas
ample.\end{prop}
\begin{proof}
Rappelons (\cite{barth}, chap. I, Lemma 17.2) que :\[
\tau_{*}\mathcal{O}_{X}=\oplus_{i=0}^{i=m-1}\mathcal{L}^{-i}\]
donc $H^{1}(X,\mathcal{O}_{X})\simeq\oplus_{i=0}^{i=m-1}H^{1}(X,\mathcal{L}^{-i})$.\\
Supposons $D=0$. Beauville \cite{BeauvilleH1} a montré que du
fait que $S$ ne possède pas de fibrations, il n'existe qu'un nombre
fini de fibrés inversibles $\mathcal{L}\in Pic^{0}(S)$ tels que $H^{1}(S,\mathcal{L})=0$.
Ainsi il existe une infinité de données de revêtements $\Lambda=(0,\mathcal{L},m)$
tels que $H^{1}(S,\mathcal{L}^{-i})=0$ pour $i=1,..,m-1$. \\
Si $D$ est ample, les groupes $H^{1}(S,\mathcal{L}^{-i})$ sont
également nuls par le théorème d'annulation de Mumford.\\
Dans les deux cas, le morphisme injectif $\tau^{*}:H^{0}(\Omega_{S})\rightarrow H^{0}(\Omega_{X})$
est un isomorphisme. Considérons le diagramme suivant:\[
\begin{array}{ccccccc}
0 & \rightarrow & \mathcal{N}\rightarrow & T_{X} & \stackrel{d\tau}{\rightarrow} & \tau^{*}T_{S}\\
 &  &  & \downarrow &  & \downarrow\\
 &  &  & H^{0}(\Omega_{X})^{*} & \rightarrow & \tau^{*}H^{0}(\Omega_{S})^{*}\end{array}\]
où les flèches verticales sont les morphismes $pr\circ d\vartheta$
définis au paragraphe \ref{sub:Les-d=0000E9finitions-de} et où $\mathcal{N}$
est le noyau de la différentielle de $\tau$. Le support de $\mathcal{N}$
est le lieu de ramification $D'$ et la suite :\[
0\rightarrow\mathcal{N}\rightarrow T_{X}\rightarrow H^{0}(\Omega_{X})^{*}\]
est exacte. \\
Si $D=0$, le corollaire \ref{diff=0000E9rentielle de l'apli. cot.}
implique que l'image d'un point $x$ de $X$ par le morphisme $\mathcal{G}_{X}$
est la même que l'image du point $\tau(x)$ par $\mathcal{G}_{S}$.
Ce morphisme $\mathcal{G}_{X}$ se factorise donc par $\tau:X\rightarrow S$
qui est de degré $m$. Les applications cotangentes de $X$ et $S$
ont donc la même image.\\
Si $D$ est ample, alors $\Omega_{X}$ n'est pas engendré sur $D'$.
Le théorème 1 de Spurr \cite{Spurr2} montre que $\Omega_{X}$ n'est
pas ample.\end{proof}
\begin{rem}
\label{rem:1)-La-d=0000E9monstration}On ne connait pas d'exemples
de surfaces ne possédant qu'un nombre fini de courbes non-amples (voir
définition supra) $D$ telles que $D^{2}>0$. La démonstration du
théorème 1 de \cite{Spurr2} montre que si $\Omega_{S}$ est ample,
la courbe $D'$ est l'unique obstruction à l'amplitude de $\Omega_{X}$.
\end{rem}
La démonstration de la proposition \ref{pro:a)-revetement etales}
montre que les seuls morphismes $X\rightarrow S$ entre deux surfaces
$X,S$ vérifiant l'hypothèse \ref{conditions} et de même irrégularité
sont étales. Les application de Gauss des deux surfaces ont alors
la même image et que les applications cotangentes des deux surfaces
également.\\
La question se pose de savoir si l'image du morphisme de Gauss
est toujours birationelle à une surface de même irrégularité. L'exemple
suivant montre que cela n'est pas le cas:

Soit $S$ une surface vérifiant les hypothèses \ref{conditions}.
Supposons que $S$ possède une involution $\iota:S\rightarrow S$
telle que le morphisme quotient $\eta:S\rightarrow X:=S/\iota$ soit
étale et telle que la surface quotient $X$ soit régulière.
\begin{prop}
Le morphisme de Gauss de $S$ se factorise par la surface régulière
$X$ et est de degré au moins $2$.\end{prop}
\begin{proof}
Le revêtement étale $\eta$ vérifie $\eta_{*}\mathcal{O}_{S}=\mathcal{\mathcal{O}}_{X}\oplus\mathcal{L}$
où $\mathcal{L}$ est un fibré inversible tel que $\mathcal{L}^{\otimes2}=\mathcal{O}_{X}$.
La codifferentielle:\[
\eta^{*}\Omega_{X}\rightarrow\Omega_{S}\]
est un isomorphisme et $\eta_{*}\eta^{*}\Omega_{X}=\Omega_{X}\otimes(\mathcal{\mathcal{O}}_{X}\oplus\mathcal{L})$
donc:\[
H^{0}(S,\Omega_{S})=H^{0}(S,\eta^{*}\Omega_{X})\simeq H^{0}(X,\Omega_{X}\otimes(\mathcal{\mathcal{O}}_{X}\oplus\mathcal{L}))=H^{0}(X,\Omega_{X}\otimes\mathcal{L}).\]
Notons $\delta_{X}$ le sous-systeme linéaire du système canonique
de $X$ défini par l'image du morphisme :\[
\wedge^{2}H^{0}(X,\Omega_{X}\otimes\mathcal{L})\rightarrow H^{0}(X,\omega_{X}).\]
Le sous-systeme linéaire $\delta_{S}$ du système canonique de $S$
défini par l'image du morphisme :\[
\wedge^{2}H^{0}(S,\Omega_{S})=\wedge^{2}H^{0}(X,\Omega_{X}\otimes\mathcal{L})\rightarrow H^{0}(S,\omega_{S})\]
vérifie: $\delta_{S}=\eta^{*}\delta_{X}$. Puisque $\Omega_{S}$ est
engendré par sections globales, le système $\delta_{S}$ est sans
points base et définit un morphisme $S\rightarrow\delta_{S}$. Ce
morphisme se factorise par $\eta$ suivi du morphisme naturel $X\rightarrow\delta_{X}\simeq\eta^{*}\delta_{S}$.
Le morphisme de Gauss \[
\mathcal{G}:S\rightarrow G(2,H^{0}(\Omega_{S})^{*})\]
se factorise donc par la surface régulière $X$ et est de degré au
moins $2$.
\end{proof}
Soit $A$ une variété abélienne de dimension $q>3$. Soit $\iota:A\rightarrow A$
une involution telle que $\iota^{*}\omega=-\omega$ pour toute $1$-forme
holomorphe $\omega$ de $A$. Une surface $S\hookrightarrow A$ lisse
d'irrégularité $q$ vérifie l'hypothèse \ref{conditions}. Supposons
de plus que $S$ soit non dégénérée dans $A$, stable par $\iota$
et ne contienne pas de points fixes de $\iota$. Le quotient de $S$
par $\iota$ est alors une surface régulière. On a ainsi caractérisé
les surfaces possédant une telle involution et dont le morphisme d'Albanese
est un plongement.
\begin{example}
On obtient une telle surface par intersection complète de $q-2$ diviseurs
$H_{1},..,H_{q-2}$ irréductibles, amples, tels que $\iota^{*}H_{i}=H_{i}$
et ne contenant pas de points fixes de l'involution $\iota$. Le théorème
de Lefschetz (\cite{Larzarsfeld} Remark 3.1.32) relatif à l'intersection
de diviseurs amples implique alors que la surface est irréductible
et d'irrégularité $q$. 
\end{example}
La famille de ces surfaces n'est pas limitée.\\
On ignore si une surface intersection complète de diviseurs amples
génériques de $A$ possède un morphisme de Gauss de degré $1$. Si
tel est le cas, cela donnerait un exemple d'une famille de surfaces
pour lequel le degré du morphisme de Gauss varie.

\subsubsection{Degré du cycle $\psi_{*}\pi^{*}C$.\label{subsection:Degr=0000E9-du-cycle}\protect \\
}

La proposition suivante met en rapport les propriétés numériques d'une
courbe $C\hookrightarrow S$ avec les propriétés numériques du cycle
$\psi_{*}\pi^{*}C$ : 
\begin{prop}
\label{intersection KC}Soit $C$ une courbe contenue dans $S$ ;
le degré de $\psi_{*}\pi^{*}C$ est égal à $KC$ où $K$ est un diviseur
canonique.\end{prop}
\begin{proof}
Soit $H$ une section hyperplane et $h=\psi^{*}H$. Par définition,
le degré de $\psi_{*}\pi^{*}C$ est le degré de l'intersection de
$H^{2}$ et de $\psi_{*}\pi^{*}C$. Le morphisme $\psi$ étant propre,
on peut utiliser la formule de projection :\[
\psi_{*}(h^{2}.\pi^{*}C)=H^{2}.\psi_{*}\pi^{*}C.\]
Les degrés de $h^{2}.\pi^{*}C$ et de $\psi_{*}(h^{2}.\pi^{*}C)$
étant égaux, il reste à montrer que le degré de $h^{2}.\pi^{*}C$
est égal à $KC$. Les classes de Chern $c_{1}(\Omega_{S}),c_{2}(\Omega_{S})$
du fibré $\Omega_{S}$ vérifient : $h^{2}=h.\pi^{*}c_{1}(\Omega_{S})-\pi^{*}c_{2}(\Omega_{S})$
dans l'anneau de Chow de $\mathbb{P}(T_{S})$. D'où :\[
h^{2}.\pi^{*}C=h.\pi^{*}c_{1}(\Omega_{S}).\pi^{*}C-\pi^{*}c_{2}(\Omega_{S}).\pi^{*}C=h.\pi^{*}c_{1}(\Omega_{S}).\pi^{*}C,\]
or le degré de $c_{1}(\Omega_{S}).C$ est $KC$.
\end{proof}

\section{Etude des courbes non-amples.}

\subsection{Critère de contraction, décomposition du fibré cotangent.\label{section2.1}}

Soit $S$ une surface vérifiant les hypothèses \ref{conditions} et
soit $C$ une courbe de $S$ i.e. un schéma de dimension $1$ réduit
et irréductible. Notons $K$ un diviseur canonique de $S$.
\begin{lem}
\label{fibretrivial}\label{fibretrivial2}Une courbe $C\hookrightarrow S$
possède une section $t:C\rightarrow\mathbb{P}(T_{S})$ contractée
en un point par l'application cotangente si et seulement si il existe
un morphisme surjectif : $\Omega_{S}\otimes\mathcal{O}_{C}\stackrel{q}{\rightarrow}\mathcal{O}_{C}\rightarrow0.$\\
En ce cas le noyau du morphisme $q$ est isomorphe à $\mathcal{O}_{C}(K)$
et la suite exacte : \[
0\rightarrow\mathcal{O}_{C}(K)\rightarrow\Omega_{S}\otimes\mathcal{O}_{C}\stackrel{q}{\rightarrow}\mathcal{O}_{C}\rightarrow0\]
 est scindée. \end{lem}
\begin{proof}
Soit $t:C\rightarrow\mathbb{P}(T_{S})$ une section et soit $\Omega_{S}\otimes\mathcal{O}_{C}\rightarrow\mathcal{L}\rightarrow0$
le quotient correpondant à la section $t$, où le fibré inversible
$\mathcal{L}$ vérifie : $t^{*}(\mathcal{O}_{\mathbb{P}(T_{S})}(1))=\mathcal{L}$
et $(\psi\circ t)^{*}(\mathcal{O}_{\mathbb{P}^{q-1}}(1))=\mathcal{L}$.
\\
Le morphisme $\psi\circ t$ est constant si et seulement si $\mathcal{L}$
est trivial. Supposons qu'un tel quotient trivial $q$ existe. Puisque
le fibré $\Omega_{S}\otimes\mathcal{O}_{C}$ est engendré par ses
sections globales, existe une section $t\in H^{o}(C,\Omega_{S}\otimes\mathcal{O}_{C})$
telle que $q(t)$ soit une section non nulle de $\mathcal{O}_{C}$
: cela fournit une section du quotient $q$.\end{proof}
\begin{defn}
Une courbe $C\hookrightarrow S$ est dite \textbf{non-ample} si la
restriction du fibré cotangent à $C$ possède un quotient isomorphe
à $\mathcal{O}_{C}$. 
\end{defn}

\subsection{Classification et lissité des courbes non-amples.\label{parag entomologie}}

Soit $C$ une courbe réduite irréductible contenue dans la surface
$S$ et $\Omega_{C}$ le faisceau des différentielles. La suite naturelle
suivante : \[
0\rightarrow\mathcal{O}_{C}(-C)\rightarrow\Omega_{S}\otimes\mathcal{O}_{C}\rightarrow\Omega_{C}\rightarrow0\]
est exacte (cf. \cite{hart} proposition 8.12) ; on l'appellera la
\textbf{suite cotangente} de la courbe $C$.\\
Le théorème suivant classifie partiellement les courbes non-amples
$C$ suivant la valeur de l'intersection $C^{2}$ :
\begin{thm}
\label{entomologie}Soit $C\hookrightarrow S$ une courbe réduite
et irréductible.\\
1) La courbe $C$ est une courbe non-ample et vérifie $C^{2}<0$
si et seulement si $C$ est lisse et de genre $1$. En ce cas la suite
cotangente est scindée.\\
2) La courbe $C$ est une courbe non-ample et vérifie $C^{2}=0$
si et seulement si $C$ est lisse, de genre $>1$ et la suite cotangente
est scindée. En ce cas le fibré normal $\mathcal{O}_{C}(C)$ est trivial.
\end{thm}
La courbe $C$ est lisse si et seulement si le faisceau $\Omega_{C}$
est localement libre de rang $1$ (\cite{hart} théorème 8.17, Chap.
II). Pour démontrer le théorème \ref{entomologie}, nous aurons besoin
d'un résultat de Lipman qui est la version duale de ce critère de
lissité.
\begin{thm}
\label{lipman2}(Lipman \cite{lipman} theorem 1) Une courbe $C$
définie sur un corps de caractéristique nulle est lisse si et seulement
si $T_{C}:=Hom_{\mathcal{O}_{C}}(\Omega_{C},\mathcal{O}_{C})$ est
un fibré inversible.
\end{thm}
Nous utiliserons également le lemme suivant (\cite{barth}, lemme
$12.2$ , Chap. II) :
\begin{lem}
\label{bricoles} Notons $\mathcal{L}$ un fibré inversible de degré
négatif ou nul sur une courbe $C$. Le fibré $\mathcal{L}$ est isomorphe
au fibré trivial si et seulement si l'espace $H^{o}(C,\mathcal{L})$
est non nul.
\end{lem}
Montrons le théorème \ref{entomologie} :

Soit $C$ une courbe contenue dans la surface. La suite duale de la
suite exacte : \[
\mathcal{O}_{C}(-C)\rightarrow\Omega_{S}\otimes\mathcal{O}_{C}\rightarrow\Omega_{C}\rightarrow0\]
 est la suite exacte : \[
0\rightarrow T_{C}\stackrel{i}{\rightarrow}T_{S|C}\stackrel{q'}{\rightarrow}\mathcal{O}_{C}(C).\]
Si $C$ est une courbe non-ample, alors la suite exacte : \[
0\rightarrow\mathcal{O}_{C}(K)\rightarrow\Omega_{S}\otimes\mathcal{O}_{C}\rightarrow\mathcal{O}_{C}\rightarrow0\]
est scindée (proposition \ref{fibretrivial}) ; la suite duale est
: \[
0\rightarrow\mathcal{O}_{C}\stackrel{i'}{\rightarrow}T_{S|C}\stackrel{q}{\rightarrow}\mathcal{O}_{C}(-K)\rightarrow0.\]
Considérons le diagramme suivant à lignes et colonnes exactes : \[
\begin{array}{cccccc}
 &  &  & \begin{array}{c}
0\\
\downarrow\end{array}\\
 &  &  & \mathcal{O}_{C} & t'\\
 &  &  & \downarrow i' & \searrow\\
0\rightarrow & T_{C} & \stackrel{i}{\rightarrow} & T_{S|C} & \stackrel{q'}{\rightarrow} & \mathcal{O}_{C}(C)\\
 &  & \searrow & \downarrow q\\
 &  & t & \mathcal{O}_{C}(-K)\\
 &  &  & \begin{array}{c}
\downarrow\\
0\end{array}\end{array}\]
où $t$ et $t'$ rendent le diagramme commutatif. 
\begin{notation}
Lorsque dans la suite de ce paragraphe $C$ dénote une courbe non-ample
de $S$, les notations $i,i',t...$ renvoient à ce diagramme.
\end{notation}
Le lemme suivant montre la première affirmation du théorème \ref{entomologie}
:
\begin{lem}
\label{lemme des 4 =0000E9quivalences}Soit $C$ une courbe contenue
dans $S$. Les quatres assertions suivantes sont équivalentes : \\
a) la courbe $C$ est non-ample et le morphisme $t$ est nul.\\
b) la courbe $C$ est non-ample et le morphisme $t'$ est nul.\\
c) la courbe $C$ est une courbe lisse de genre $1$.\\
d) la courbe $C$ est non-ample et $C^{2}<0$.\\
Si une des assertions est vérifiée, alors la suite cotangente de
$C$ est scindée.\end{lem}
\begin{proof}
Soit $C\hookrightarrow S$ une courbe. Supposons que le point a) soit
vérifié i.e. $C$ est non-ample et le morphisme : \[
t=q\circ i:T_{C}\rightarrow\mathcal{O}_{C}(-K)\]
 est nul. Montrons que b) est vérifié. Puisque $t=q\circ i$ est nul,
par propriété du noyau de $q$, il existe un morphisme $j:T_{C}\rightarrow\mathcal{O}_{C}$
tel que le diagramme suivant :\[
\begin{array}{ccc}
 &  & \mathcal{O}_{C}\\
 & j\nearrow & \downarrow i'\\
T_{C} & \stackrel{i}{\rightarrow} & T_{S|C}\end{array}\]
commute et tel que $j$ soit injectif (car $i$ est injectif). Puisque
le morphisme $t'$ se factorise par $q'$, $t'$ est nul sur le sous-$\mathcal{O}_{C}$-module
$j(T_{C})$ de $\mathcal{O}_{C}$. Le morphisme $t'$ se factorise
donc par le quotient de $\mathcal{O}_{C}$ par $j(T_{C})$. Ce quotient
est de torsion et puisque $\mathcal{O}_{C}(C)$ est un fibré inversible,
le morphisme : \[
t':\mathcal{O}_{C}\rightarrow\mathcal{O}_{C}(C)\]
 est nul. Nous avons donc montré que a) implique b).\\
Réciproquement montrons que b) implique a). La situation est symétrique
de la précédente implication : \\
Soit $C$ une courbe non-ample telle que le morphisme $t'=q'\circ i'$
soit nul. Par propriété du noyau de $q'$, l'injection $i':\mathcal{O}_{C}\rightarrow T_{S|C}$
se factorise par un morphisme $j':\mathcal{O}_{C}\rightarrow T_{C}$
tel que le diagramme suivant :\[
\begin{array}{ccc}
 &  & \mathcal{O}_{C}\\
 & j'\swarrow & \downarrow i'\\
T_{C} & \stackrel{i}{\rightarrow} & T_{S|C}\end{array}\]
commute. Le morphisme $j'$ est de plus injectif car $i'$ est injectif.
Puisque le morphisme $t$ se factorise par $q$, $t$ est nul sur
le sous-$\mathcal{O}_{C}$-module $j'(\mathcal{O}_{C})$ de $T_{C}$.
Le morphisme $t$ se factorise donc par le quotient de $T_{C}$ par
$j'(\mathcal{O}_{C})$. Ce quotient est de torsion et puisque $\mathcal{O}_{C}(-K)$
est un fibré inversible, le morphisme : \[
t:T_{C}\rightarrow\mathcal{O}_{C}(-K)\]
est nul. Nous avons donc montré que b) implique a).\\
Il y a donc équivalence entre les points a) et b) et de plus, si
$C$ est une courbe non-ample vérifiant a) et b), nous avons construit
deux morphismes injectifs $j:T_{C}\rightarrow\mathcal{O}_{C}$ et
$j':\mathcal{O}_{C}\rightarrow T_{C}$. Le composé $j\circ j':\mathcal{O}_{C}\rightarrow\mathcal{O}_{C}$
est un morphisme injectif, c'est donc un isomorphisme et on en déduit
que $j:T_{C}\rightarrow\mathcal{O}_{C}$ est un morphisme surjectif.
Le faisceau $T_{C}$ est donc isomorphe à $\mathcal{O}_{C}$. Le théorème
de Lipman \ref{lipman2} permet alors de conclure que la courbe $C$
est lisse. \\
En ce cas $\Omega_{C}$ est un fibré inversible, donc : \[
\hom_{\mathcal{O}_{C}}(T_{C},\mathcal{O}_{C})=\Omega_{C},\]
 et puisque $T_{C}$ est trivial, le fibré $\Omega_{C}$ est trivial.
Ainsi $C$ est lisse de genre $1$ et nous avons montré que si le
point a) ou b) est vérifié, alors c) est vérifié.\\
Soit $C\hookrightarrow S$ une courbe vérifiant c) i.e. $C$ est
lisse de genre $1$. Le quotient naturel $\Omega_{S|C}\rightarrow\Omega_{C}$
est un quotient trivial et surjectif donc $C$ est une courbe non-ample.
De plus, par adjonction : $C^{2}+KC=0$ et puisque les hypothèses
sur $S$ entrainement que le diviseur canonique $K$ est ample, on
en déduit que nécessairement $C^{2}$ est strictement négatif. Ceci
montre que c) entraîne d). (Remarquons de plus que la proposition
\ref{fibretrivial} montre que la suite cotangente est scindée).\\
Soit $C$ une courbe vérifiant l'hypothèse d) i.e. $C$ est non-ample
et $C^{2}<0$. En ce cas, le morphisme $t':\mathcal{O}_{C}\rightarrow\mathcal{O}_{C}(C)$
est une section du fibré $\mathcal{O}_{C}(C)$. Mais ce fibré est
de degré $C^{2}$ strictement négatif et le lemme \ref{bricoles}
montre que $t'=0$. Nous avons montré que d) entraîne b).\\
Les quatres assertions a), b), c), d) sont donc équivalentes.
\end{proof}
Montrons maintenant la seconde affirmation du théorème \ref{entomologie}.

- Soit $C\hookrightarrow S$ une courbe non-ample vérifiant $C^{2}=0$,
montrons que $C$ est lisse de genre $>1$ et que la suite cotangente
est scindée. Le morphisme : \[
t'=q'\circ i':\mathcal{O}_{C}\rightarrow\mathcal{O}_{C}(C)\]
 est nécessairement non nul car nous avons montré au lemme \ref{lemme des 4 =0000E9quivalences}
que $t'=0$ entraîne $C^{2}<0$. Ce morphisme $t'$ peut être considéré
comme une section non nulle de l'espace $H^{o}(C,\mathcal{O}_{C}(C))$.
Puisque le fibré $\mathcal{O}_{C}(C)$ est de degré $C^{2}=0$, le
lemme \ref{bricoles} montre que $\mathcal{O}_{C}(C)$ est trivial.
Le morphisme : \[
t'=q'\circ i':\mathcal{O}_{C}\rightarrow\mathcal{O}_{C}(C)\simeq\mathcal{O}_{C}\]
 est alors un isomorphisme ; le morphisme : \[
q':T_{S|C}\rightarrow\mathcal{O}_{C}(C)=\mathcal{O}_{C}\]
 est donc surjectif et son noyau $T_{C}$ est un fibré inversible.
Nous pouvons donc appliquer le théorème de Lipman \ref{lipman2} et
conclure que si $C$ est une courbe non-ample telle que $C^{2}=0$,
alors la courbe $C$ est lisse. De plus, la suite suivante : \[
0\rightarrow T_{C}\rightarrow T_{S|C}\rightarrow\mathcal{O}_{C}(C)=\mathcal{O}_{C}\rightarrow0\]
est exacte et scindée par un multiple du morphisme $i':\mathcal{O}_{C}\rightarrow T_{S|C}$
(car $t'=q'\circ i'$ est un isomorphisme). La suite cotangente qui
est donc scindée.

- Réciproquement, soit $C$ une courbe lisse contenue dans la surface
telle que $C^{2}=0$ et telle que la suite cotangente soit scindée
; montrons que $C$ est non-ample. Puisque la suite cotangente est
scindée, le fibré $\Omega_{S}\otimes\mathcal{O}_{C}$ est isomorphe
à $\mathcal{O}_{C}(-C)\oplus\Omega_{C}$. Mais $\Omega_{S}\otimes\mathcal{O}_{C}$
est engendré par restriction des sections globales de $\Omega_{S}$,
donc l'espace $H^{o}(C,\mathcal{O}_{C}(-C))$ est non nul. Par le
lemme \ref{bricoles}, $\mathcal{O}_{C}(-C)$ est trivial. Ainsi le
fibré $\Omega_{S}\otimes\mathcal{O}_{C}\simeq\mathcal{O}_{C}\oplus\Omega_{C}$
admet un quotient trivial et $C$ est non-ample. 

Soit $C$ lisse non-ample vérifiant $C^{2}>0$. La suite cotangente
: \[
0\rightarrow\mathcal{O}_{C}(-C)\rightarrow\Omega_{S}\otimes\mathcal{O}_{C}\rightarrow\Omega_{C}\rightarrow0\]
 ne peut être scindée car $\Omega_{S}\otimes\mathcal{O}_{C}$ est
engendré par ses sections globales mais le fibré $\mathcal{O}_{C}(-C)$
est de degré $-C^{2}<0$. 

Ceci achève la démonstration du théorème \ref{entomologie}. $\Box$

\subsection{Exemple des fibrations et géométrie de l'image par $\psi$ de $\pi^{*}C$.\label{sub:G=0000E9om=0000E9trie-de-l'image}\protect \\
}

Soit $S$ une surface vérifiant les hypothèses \ref{conditions} et
telle qu'il existe un morphisme $f:S\rightarrow B$ surjectif à fibres
connexes dans une courbe lisse $B$. En \cite{Martin} p. 50, Martin-Deschamps
montre qu'une fibre lisse $f^{*}b$ ($b$ point de $B$) est non-ample
si et seulement si $b$ est un zéro du morphisme de Kodaira-Spencer
$\delta$ associé à la fibration (pour la définition de $\delta$
voir \cite{Szpiro}, exposé III) . La fibration est dite isotriviale
si le morphisme de Kodaira-Spencer est nul. Ainsi:
\begin{cor}
\label{cor:Quand-la-fibration}Si la fibration est isotriviale, alors
la surface possède une infinité de courbes non-amples.
\end{cor}
On dispose ainsi d'exemples de surfaces ayant un nombre infini de
courbes non-amples telles que $C^{2}=0$. 

Soit $C$ une courbe contenue $S$. Notons $T$ la surface image de
$\pi^{*}C$ par $\psi$, notons $\mathcal{K}$ le noyau du morphisme
de restriction : \[
H^{o}(\Omega_{S})\rightarrow H^{o}(C,\Omega_{S}\otimes\mathcal{O}_{C})\]
et $k$ la dimension de $\mathcal{K}$. Notons de plus $V$ l'espace
quotient de $H^{0}(\Omega_{S})$ par $\mathcal{K}$ et $V^{*}$ son
dual. L'espace projectif $\mathbb{P}(V^{*})$ est naturellement plongé
dans $\mathbb{P}(H^{0}(\Omega_{S})^{*})$. 
\begin{prop}
\label{prop la surface psi(pi(c))}A) L'enveloppe linéaire de la surface
$T$ dans $\mathbb{P}^{q-1}$ est l'espace projectif $\mathbb{P}(V^{*})\hookrightarrow\mathbb{P}^{q-1}$
de dimension $q-k-1$.\\
B) Si $\mathbb{P}(V^{*})\not=\mathbb{P}^{q-1}$, alors $C$ vérifie
: $C^{2}\leq0$. Si de plus $C^{2}=0$, alors un multiple de $C$
est une fibre d'une fibration de $S$ dans une courbe de genre $b\geq1$.\\
C) De plus, si $C$ est une courbe non-ample, alors $b=k+1$.\\
D) Réciproquement, si $C$ est non-ample et fibre d'une fibration
de $S$ dans une courbe de genre $b\geq1$, alors $b=k+1$.
\end{prop}
Cette proposition entraine qu'une courbe non-ample $C$ telle que
$C^{2}>0$ crée une singularité sur l'image de l'application cotangente
pourvu que $C^{2}>0$ (voir également corollaire \ref{cor:sommet du cone lisse}).

Soit $i:C\hookrightarrow Z$ une courbe sur une surface lisse $Z$.
On note $\Omega_{C}$ le faisceau des différentielles de $C$. Pour
démontrer la proposition, on utilisera principalement le lemme suivant
dû à Spurr (\cite{Spurr} Theorem 1):
\begin{lem}
\label{lem:Spurr1}Soit $\omega$ une $1$-forme holomorphe telle
que $i^{*}\omega=0\in H^{0}(C,\Omega_{C})$. Alors la courbe vérifie
: $C^{2}\leq0$. Si de plus $C^{2}=0$, alors un multiple de $C$
est la fibre d'une fibration $f:Z\rightarrow B$ dans une courbe $B$
lisse de genre $b\geq1$.
\end{lem}
Ce lemme a la conséquence directe suivante :
\begin{lem}
\label{lem:Spurr2}Soit $f:Z\rightarrow B$ une fibration dans une
courbe $B$ et soit $i:C\hookrightarrow Z$ une courbe telle que $mC$
soit une fibre de $f$ (pour un certain $m>0$). Une $1$-forme holomorphe
$\tau$ de $Z$ est un élément de $f^{*}H^{0}(B,\Omega_{B})$ si et
seulement si la restriction $i^{*}\tau\in H^{0}(C,\Omega_{C})$ est
nulle.\end{lem}
\begin{proof}
(De la proposition \ref{prop la surface psi(pi(c))}). Posons $Y=\pi^{-1}(C)$
et $j:Y\hookrightarrow\mathbb{P}(T_{S})$ le morphisme d'inclusion.
Le morphisme $\psi\circ j:Y\rightarrow\mathbb{P}^{q-1}$ est obtenu
par le quotient : \[
H^{o}(\Omega_{S})\otimes\mathcal{O}_{Y}\rightarrow j^{*}\mathcal{O}_{\mathbb{P}(T_{S})}(1)\rightarrow0\]
qui se factorise comme suit :\[
H^{o}(\Omega_{S})\otimes\mathcal{O}_{Y}\rightarrow H^{o}(C,\Omega_{S}\otimes\mathcal{O}_{C})\otimes\mathcal{O}_{Y}\rightarrow({j\circ\pi)}^{*}(\Omega_{S}\otimes\mathcal{O}_{C})\rightarrow j^{*}\mathcal{O}_{\mathbb{P}(T_{S})}(1)\rightarrow0.\]
L'image de la surface $Y$ par $\psi$ est donc contenue et non-dégénérée
dans $\mathbb{P}(V^{*})\hookrightarrow\mathbb{P}^{q-1}$.\\
Si ${\mathbb{P}(V}^{*})$ est strictement contenu dans $\mathbb{P}^{q-1}$,
alors il existe une $1$-forme dont la restriction à $C$ est nulle.
Par le lemme \ref{lem:Spurr1}, on a alors : $C^{2}\leq0$ et de plus,
si $C^{2}=0$, alors $C$ est une fibre d'une fibration $f:S\rightarrow B$
dans une courbe de genre $b\geq1$.

Considérons maintenant $C\hookrightarrow S$ une courbe non-ample.
Supposons qu'un multiple $mC$ de $C$ soit la fibre d'une fibration
$f:S\rightarrow B$. Par le théorème \ref{entomologie}, la courbe
$C$ est lisse, le fibré normal est trivial et: \[
\Omega_{S}\otimes\mathcal{O}_{C}=\mathcal{O}_{C}\oplus\mathcal{O}_{C}(K)\]
où $K$ est un diviseur canonique de $S$. On identifie $\mathcal{O}_{C}(K)$
avec le fibré canonique de $C$. Si $mC=f^{*}p$ (où $p$ est un point
de $B$), le noyau $\mathcal{K}$ du morphisme de restriction: \[
H^{o}(\Omega_{S})\rightarrow H^{o}(C,\Omega_{S}\otimes\mathcal{O}_{C})=H^{0}(C,\mathcal{O}_{C})\oplus H^{0}(C,\Omega_{C})\]
est l'image réciproque par $f$ de l'hyperplan de $H^{0}(B,\Omega_{B})$
formé des formes nulles en $p$, ainsi $b=k+1$.
\end{proof}
Soit $C$ une courbe non-ample de $S$. Soit $T$ le cône image de
$\pi^{*}C$ par $\psi$. Si l'image de l'application cotangente est
lisse au sommet $t$ du cône $T$, alors les droites passant par $t$
sont contenues dans l'espace projectif tangent à $t$. Ainsi la surface
$T$ est contenue dans un sous-espace projectif de dimension $3$
de $\mathbb{P}^{q-1}$. Si $q\geq5$, la proposition \ref{prop la surface psi(pi(c))}
et le théorème \ref{entomologie} impliquent le résultat suivant:
\begin{cor}
\label{cor:sommet du cone lisse}Sous les hypothèses précédentes,
la courbe $C$ vérifie l'une des deux propriétés suivantes:\\
a) ou bien $C^{2}<0$ et $C$ est une courbe elliptique,\\
b) ou bien $C^{2}=0$ et un multiple de $C$ est la fibre d'une
fibration dans une courbe de genre $b$ et $q-3\leq b\leq q-2$.
\end{cor}
Les courbes de genre petit ont des conséquences particulières sur
l'image de l'application cotangente:
\begin{cor}
Soit $C\hookrightarrow S$ une courbe lisse de genre $2$. La courbe
$C$ vérifie : $C^{2}\leq0$. Si $C^{2}=0$, alors $C$ est une courbe
non-ample et il existe un entier $n$ tel que $nC$ soit la fibre
d'une fibration $f:S\rightarrow B$ dans une courbe de genre $b=q-3$.
\\
Si de plus $n=1$, alors la fibration est isotriviale à fibres
de genre $2$. \\
Il existe alors une droite $L\hookrightarrow\mathbb{P}^{q-1}$
et une courbe $D\hookrightarrow\mathbb{P}^{q-1}$ disjointes et telles
que l'image de l'application cotangente soit balayée par les plans
passant par un point de $D$ et contenant $L$.
\end{cor}
La dernière assertion est un cas particulier du corollaire \ref{pro:Soit-s=00003D(p_{1},p_{2})-un}
démontré plus loin.
\begin{proof}
Notons $K$ un diviseur canonique. La restriction à $C$ du morphisme
de Gauss $\mathcal{G}:S\rightarrow G(2,q)$ suivis du plongement de
Plücker de $G(2,q)$ est donné par le fibré $\mathcal{O}_{C}(K)$.
Puisque $\mathcal{G}$ est fini, on a $KC\geq2$ et si $KC=2$, alors
l'image de $C$ par $\mathcal{G}$ est une droite. Une droite de la
grassmanienne correspond à l'ensemble des droites passant par un point
et contenues dans un plan de $\mathbb{P}^{q-1}$. Ainsi $C$ est non-ample.\\
L'espace $H^{0}(C,\Omega_{S}\otimes\mathcal{O}_{C})\simeq H^{0}(C,\mathcal{O}_{C}(K)\oplus\mathcal{O}_{C})$
est de dimension $3$, la base de la fibration est donc de genre $q-3$.\\
Si $n=1$, les fibres lisses sont de genre $2$, donc non-amples,
et la fibration est isotriviale. 
\end{proof}

\section{Surfaces contenant une infinité de courbes non-amples.}

\subsection{Caractérisation de l'image de l'application cotangente.}

Rappelons que nous avons noté $\Delta$ le lieu des points exceptionnels
i.e. des points $p$ de $\mathbb{P}^{q-1}$ tels que la fibre $\psi^{-1}(p)$
soit de dimension $1$. La dimension de $\Delta$ est inférieure ou
égale à $1$.\\
Supposons que le fermé $\Delta$ contienne une composante irréductible
$B$ de dimension $1$. En ce cas son image inverse $\psi^{-1}(B)$
est de dimension $2$. Soit $S'$ une composante irréductible de dimension
$2$ de $\psi^{-1}(B)$ munie de la structure réduite. Notons $\pi'$
et $\psi'$ les restrictions à $S'$ du morphisme de projection $\pi$
et de l'application cotangente :\[
\begin{array}{ccc}
S' & \stackrel{\psi'}{\rightarrow} & B.\\
\pi'\downarrow\\
S\end{array}\]
Rappelons qu'une droite de $\mathbb{P}^{q-1}$ est dite sécante d'une
courbe $X\hookrightarrow\mathbb{P}^{q-1}$ si elle coupe $X$ en au
moins deux points. Si $Y\hookrightarrow\mathbb{P}^{q-1}$ est une
seconde courbe, une droite est dite sécante de $X$ et de $Y$ si
elle passe par $X$ et $Y$.

Le théorème suivant caractérise les surfaces qui contiennent une infinité
de courbes non-amples :
\begin{thm}
\label{thitalien}Le morphisme $\pi'$ est surjectif. Si $s$ est
un point générique de $S$ alors la droite $L_{s}$ coupe $B$ en
$d_{o}\in\{0,1,2\}$ points où $d_{o}$ est le degré de $\pi'$. \\
1) Si $d_{o}=1$ et si $B$ n'est pas une droite de $\mathbb{P}^{q-1}$,
alors le morphisme $\pi':S'\rightarrow S$ est un isomorphisme. La
surface $S$ possède une fibration isotriviale.\\
Le lieu des points exceptionnels $\Delta$ est formé de deux courbes
lisses irréductibles ${B=B}_{1},B_{2}$ qui sont images de morphismes
canoniques de deux courbes.\\
L'image de l'application cotangente est la variété des sécantes
de ces courbes et l'image du morphisme de Gauss est la surface $B_{1}\times B{}_{2}$.\\
2) Si $d_{o}=2$, alors l'image de l'application cotangente est
la variété des sécantes de $B$. Une courbe non-ample $C\hookrightarrow S$
vérifie : $C^{2}=\deg\mathcal{G}>0$ où $\deg\mathcal{G}$ est le
degré du morphisme de Gauss. \\
L'image du morphisme de Gauss est birationnelle à $B^{(2)}$.
\end{thm}
Commençons par le lemme suivant :
\begin{lem}
\label{horizon}Le morphisme $\pi':S'\rightarrow S$ est surjectif.\end{lem}
\begin{proof}
Si $S'$ est un diviseur vertical pour la projection $\pi$, alors
la courbe $B=\psi(S')$ est une droite projective et pour tout point
$s$ de $D=\pi(S')$, $L_{s}=B$. La courbe $D$ est alors contractée
en un point par le morphisme de Gauss. Cela est impossible car $\mathcal{G}$
est fini (lemme \ref{gauss est fini}). \\
Le morphisme $\pi$ étant propre, le morphisme $\pi'$ est surjectif.
\end{proof}
Le degré $d_{o}$ de $\pi'$ est égal au nombre d'intersection dans
$\mathbb{P}(T_{S})$ de $S'$ et de la fibre $\pi^{-1}s$ en un point
générique $s$ de $S$.
\begin{lem}
\label{reunions}Soit $s$ un point de la surface. La droite $L_{s}$
coupe $B$ et si $s$ est générique, alors $L_{s}$ coupe la courbe
$B$ en $d_{o}$ points. \end{lem}
\begin{proof}
Par le lemme \ref{horizon}, la fibre $\pi^{-1}(s)$ coupe $S'$ dans
$\mathbb{P}(T_{S})$ donc l'image de $\pi^{-1}(s)$ par $\psi$ coupe
l'image de $S'$ par $\psi$, c'est-à-dire : $L_{s}$ coupe $B$.
\\
L'application cotangente est un plongement sur $\pi^{-1}(s)$.
Cela implique que si l'intersection de $\pi^{-1}(s)$ et de $S'$
contient $d_{o}$ points distincts, alors $L_{s}$ coupe $B$ en $d_{o}$
points distincts.
\end{proof}

\subsubsection*{Démonstration de la partie 1) du théorème \ref{thitalien}.\protect \\
}

Supposons que le degré $d_{o}$ de $\pi':S'\rightarrow S$ vaut $1$.
\begin{lem}
\label{si do =00003D 1 alors S est isotriviale}Si $d_{o}=1$ et si
$B$ n'est pas une droite, alors le morphisme $\pi':S'\rightarrow S$
est un isomorphisme et la fibration de Stein associée à la fibration
$\psi\circ\pi'^{-1}:S\rightarrow B$ est isotriviale.\end{lem}
\begin{proof}
Puisque le degré $d_{o}$ de $\pi'$ est égal à $1$, le morphisme
$\pi':S'\rightarrow S$ est birationnel (cf. \cite{milne} remarque
6.21). La surface $S$ est normale. Le théorème principal de Zariski
montre que si l'application rationnelle $\pi'^{-1}$ n'est pas définie
en un point $s_{o}$ de $S$, alors l'image inverse $\pi'^{-1}s_{o}$
est une courbe de $S'\hookrightarrow\mathbb{P}(T_{S})$. \\
En ce cas, la fibre $\pi^{-1}s_{o}$ est contenue dans $S'$. Puisque
l'application cotangente est un plongement sur $\pi^{-1}s_{o}$, on
en déduit alors que $B=\psi'(S')$ est égal à $L_{s_{o}}$. Mais le
cas 1) suppose que $B$ n'est pas une droite de $\mathbb{P}^{q-1}$.\\
 Il existe donc un morphisme réciproque $t:S\rightarrow S'$ à
$\pi'$. Les composantes lisses irréductibles d'une fibre de $\psi\circ t:S\rightarrow B$
sont des courbes non-amples. Le corollaire \ref{cor:Quand-la-fibration}
permet conclure que la surface $S$ admet une fibration isotriviale.
\end{proof}
Il nous reste à caractériser l'image de l'application cotangente.
Supposons que la courbe $B$ ne soit pas une droite et que le degré
de $\pi':S'\rightarrow S$ soit égal à $1$, alors :
\begin{prop}
\label{deuxcomposantes} Le fermé $\Delta$ est formé de deux courbes
lisses irréductibles $C'_{1},C'_{2}$ qui sont images de morphismes
canoniques de courbes.\\
L'image de $\psi$ est la réunion des droites sécantes des courbes
$C_{1}'$ et $C_{2}'$.
\end{prop}
Pour montrer cette proposition, nous allons nous ramener à supposer
que $S$ est un produit de deux courbes.

Soit $S$ une surface vérifiant les hypothèses de la proposition \ref{deuxcomposantes}.
Le lemme \ref{si do =00003D 1 alors S est isotriviale} montre que
$S$ est isotriviale. Il existe en ce cas un morphisme à fibres connexes
$f:S\rightarrow C'_{2}$ dont les fibres lisses sont isomorphes entre
elles, notons $C_{1}$ une telle fibre, alors : 
\begin{lem}
\label{lemme caract=0000E9risant les surfaces isotriv}Il existe une
courbe lisse $C_{2}$ et un groupe $G$ agissant algébriquement sur
$C_{1}$ et $C_{2}$ tels que :\\
i) La surface $S$ est birationnelle au quotient $(C_{1}\times C_{2})/G$.\\
ii) La courbe $C'_{2}$ est isomorphe à $C_{2}/G$.\\
iii ) Le diagramme suivant commute : 

\[
\begin{array}{ccc}
S & -\rightarrow & (C_{1}\times C_{2})/G\\
\downarrow &  & \downarrow\\
C'_{2} & \simeq & C_{2}/G.\end{array}\]
Ici le groupe $G$ agit sur $C_{1}\times C_{2}$ composante par composante
(i.e. $\gamma.(a,b)=(\gamma a,\gamma b)$) et la flèche verticale
de droite est la projection naturelle. \end{lem}
\begin{proof}
Voir \cite{serrano}. 
\end{proof}
Soit $S$ une surface vérifiant l'hypothèse de la proposition \ref{deuxcomposantes}
et : \[
C_{1},C_{2},G,C'_{2}\]
 possédant les propriétés du lemme \ref{lemme caract=0000E9risant les surfaces isotriv}.\\
Puisque la surface $S$ ne contient pas de courbes rationelles,
la surface $(C_{1}\times C_{2})/G$ est lisse et égale à $S$.\\
La surface isotriviale $(C_{1}\times C_{2})/G$ admet deux fibrations
$g:S\rightarrow C'_{1}:=C_{1}/G$ et $f:S\rightarrow C'_{2}$, notons
$\tau:S\rightarrow C'_{1}\times C'_{2}$ le morphisme $\tau=(g,f)$.
Le lemme suivant est la proposition 2.2 de \cite{serrano} : 
\begin{lem}
L'espace des sections globales du fibré cotangent est : \[
H^{o}(\Omega_{S})=\tau^{*}(H^{o}(C'_{1},\Omega_{C'_{1}})\oplus H^{o}(C'_{2},\Omega_{C'_{2}})).\]

\end{lem}
Les surfaces $S$ et ${C'}_{1}\times C'_{2}$ ont la même irrégularité.
Par la remarque \ref{rem:1)-La-d=0000E9monstration}, on obtient : 
\begin{cor}
\label{cor:Par-l'identification}Par l'identification $H^{o}(S,\Omega_{S})\simeq H^{o}(C'_{1}\times C'_{2},\Omega_{C'_{1}\times C'_{2}})$,
les applications cotangentes de $S$ et du produit $C'_{1}\times C'_{2}$
ont la même image. 
\end{cor}
Pour terminer la démonstration de la proposition \ref{deuxcomposantes},
il nous reste à comprendre quelle est l'image de l'application cotangente
d'une surface isotriviale $S=C_{1}\times C_{2}/G$ quand $G=\{1\}$
est le groupe trivial. Puisque $S$ est de type général, ces deux
courbes sont de genres respectifs $g_{1},g_{2}$ supérieurs ou égaux
à $2$. le lemme suivant est classique: 
\begin{lem}
Soit $C_{1},C_{2}$ deux courbes de genre $g_{1}>1$ et $g_{2}>1$.
La surface $C_{1}\times C_{2}$ vérifie l'hypothèse \ref{conditions}
et ses nombres de Chern vérifient : \[
\begin{array}{c}
c_{1}^{2}[S]=8(g_{1}-1)(g_{2}-1)\\
c_{2}[S]=4(g_{1}-1)(g_{2}-1)\end{array}\]

\end{lem}
Soit donc $S=C_{1}\times C_{2}$ avec $C_{1},C_{2}$ de genre $g_{1}>1$
et $g_{2}>1$. Notons $\pi_{i}:S\rightarrow C_{i},\, i\in\{1,2\}$
les projections respectives. Le fibré cotangent vérifie : \[
\Omega_{S}=\pi_{1}^{*}\Omega_{C_{1}}\oplus\pi_{2}^{*}\Omega_{C_{2}},\]
et l'espace des sections globales $H^{o}(\Omega_{S})$ s'identifie
à $H^{o}(C,\Omega_{C_{1}})\oplus H^{o}(C_{2},\Omega_{C_{2}})$. La
surface $S=C_{1}\times C_{2}$ est d'irrégularité $g_{1}+g_{2}>3$.\\
Pour $i\in\{1,2\}$, notons $\phi_{i}:C_{i}\rightarrow\mathbb{P}^{g_{1}+g_{2}-1}$
le composé du morphisme canonique : \[
C_{i}\rightarrow\mathbb{P}(H^{o}(C_{i},\Omega_{C_{i}})^{*})\]
avec le plongement naturel : \[
\mathbb{P}(H^{o}(C_{i},\Omega_{C_{i}})^{*})\hookrightarrow\mathbb{P}(H^{o}(C_{1},\Omega_{C_{1}})^{*}\oplus H^{o}(C_{2},\Omega_{C_{2}})^{*})=\mathbb{P}^{g_{1}+g_{2}-1}.\]
Notons de plus $C'_{i}$ l'image du morphisme $\phi_{i}$. La proposition
suivante caractérise l'image de l'application cotangente de $S$ et
est une conséquence directe de la définition du morphisme de Gauss
\ref{d=0000E9finition du morphisme de gauss}:
\begin{prop}
\label{pro:Soit-s=00003D(p_{1},p_{2})-un}Soit $s=(p_{1},p_{2})$
un point de la surface $S=C_{1}\times C_{2}$. La droite $L_{s}$
passe par les points $\phi_{1}(p_{1})$ et $\phi_{2}(p_{2})$ ; les
courbes lisses $C_{1}'$ et $C_{2}'$ forment l'ensemble des points
exceptionnels.
\end{prop}
La première partie du théorème \ref{deuxcomposantes} est donc démontrée.

\subsubsection*{Démonstration de la partie 2) du théorème \ref{thitalien}.\protect \\
}

Supposons maintenant que le degré $d_{o}$ du morphisme $\pi':S'\rightarrow S$
soit supérieur ou égal à $2$.

Soit $s$ un point générique de $S$. Par le lemme \ref{reunions},
la droite $L_{s}$ passe par $d_{o}>1$ points de $B$ et il y a trois
possibilités :\\
i) Il existe un point $b_{1}$ de $B$ telle que la droite générique
$L_{s}$ passe par $b_{1}$.\\
ii) Il existe des points $b_{1},\not=b_{2}$ de $B$ telle que
la droite $L_{s}$ passe par $b_{o},b_{1}$.\\
iii) L'intersection de $L_{s}$ et de $B$ se fait en deux points
variables de $B$.\\
Le cas i) est exclu car alors $\psi$ aurait une fibre de dimension
$2$ en $b_{1}$, ce qui contredirait le lemme \ref{pi est injectif}.
Le cas ii) est exclu car l'image de l'application cotangente serait
alors une droite.\\
Ceci montre que l'image de $\psi$ est la variété développée par
les sécantes de $B$. Puisque l'image de $\psi$ est non-dégénérée,
la courbe $B$ est non-dégénérée dans $\mathbb{P}^{q-1}$. Les fibres
de $\pi'$ sont donc finies car si $S'$ contenait une fibre $\pi^{-1}s_{o}$,
alors $B$ serait une droite et l'image de l'application cotangente
serait cette droite. \\
Soient $b,b'$ deux points génériques de $B$. La courbe $C_{b}=\pi(\psi^{-1}(b))$
est une courbe non-ample et $C_{b}^{2}=C_{b}C_{b'}$. Une sécante
de $B$ générique est repérée de manière unique par les deux points
$b,b'$ de l'intersection de $L$ et $B$, et le cardinal des droites
$L_{s}$ telles que $L_{s}=L$ est égal à $\deg\mathcal{G}$ (voir
les diagrammes du paragraphe \ref{paragraphe degr=0000E9 appli cot}),
ainsi : $C_{b}^{2}=\deg\mathcal{G}>0.$

Ceci termine la démonstration du théorème \ref{thitalien} $\Box$.

\subsection{Compléments au cas d'une surface produit de deux courbes.}

Soit $C_{1},C_{2}$ deux courbes de genres respectifs $g_{1},g_{2}$
supérieur ou égaux à $2$ et soit $n\in\{0,1,2\}$ le nombre de courbes
hyperelliptiques parmis $C_{1},C_{2}$.
\begin{prop}
La surface $S=C_{1}\times C_{2}$ vérifie l'hypothèse \ref{conditions}.
Le degré de l'image $F$ de l'application cotangente de la surface
$C_{1}\times C_{2}$ est : \[
2^{2-n}(g_{1}-1)(g_{2}-1).\]
et le degré de l'application cotangente et du morphisme de Gauss est
$2^{n}$.\end{prop}
\begin{proof}
On reprend les notations de la démonstration du théorème \ref{thitalien}.
Soit $i\in\{1,2\}$, notons $k_{i}$ le degré de $C'{}_{i}\hookrightarrow\mathbb{P}^{q-1}$.
Si la courbe $C_{i}$ est hyperelliptique alors $k_{i}=g_{i}-1$,
sinon $k_{i}=2(g_{i}-1)$. \\
Soit $(\omega_{1},\omega_{2})\in H^{o}(C_{1},\Omega_{C_{1}})\times H^{o}(C_{2},\Omega_{C_{2}})$
deux formes génériques. L'intersection de $F$ avec l'hyperplan $H_{1}=\{\omega_{1}=0\}$
est une surface formée des $k_{1}$ cônes reliant $k_{1}$ points
de $C'_{1}$ aux points de $C'_{2}$. \\
L'intersection de $H_{1}F$ avec l'hyperplan $H_{2}=\{\omega_{2}=0\}$
est formée des droites reliant $k_{1}$ points de $C_{1}$ à $k_{2}$
points de $C_{2}$ et est de degré $k_{1}k_{2}$. Ainsi $\deg F=k_{1}k_{2}$.\\
Pour le degré de l'application cotangente, on utilise la proposition
\ref{degre de psi}.
\end{proof}

\subsection{Exemple du produit symétrique.\label{paragraphe exemple du produit sym=0000E9trique}}

Soit $C$ une courbe de genre $g>3$. L'involution $\tau:(P_{1},P_{2})\rightarrow(P_{2},P_{1})$
agit sur la surface $C\times C$. On notera $S=C^{(2)}$ la surface
quotient et \[
\begin{array}{ccc}
\eta:C\times C & \rightarrow & C^{(2)}\\
(P_{1},P_{2}) & \rightarrow & P_{1}+P_{2}\end{array}\]
le morphisme quotient. Il existe un isomorphisme naturel :\[
H^{o}(C,\Omega_{C})\rightarrow H^{o}(S,\Omega_{S})\]
appelé morphisme trace \cite{Griffiths} qui permet d'identifier les
deux espaces $H^{o}(\Omega_{S})$ et $H^{o}(C,\Omega_{C})$. Puisque
la courbe $C$ n'est pas hyperelliptique, le morphisme canonique :
\[
\phi:C\rightarrow\mathbb{P}^{g-1}=\mathbb{P}(H^{o}(\Omega_{S})^{*})\]
est un plongement. 
\begin{prop}
\label{proposition produit sym=0000E9trique}Soit $C$ une courbe
non-hyperelliptique de genre $g>3$. \\
La surface $S=C^{(2)}$ vérifie l'hypothèse \ref{conditions}.
\\
Pour tout point $P$ de $C$, la courbe $P+C\hookrightarrow C^{(2)}$
est une courbe non-ample et vérifie $(P+C)^{2}=1$.\\
L'image de l'application cotangente est la variété des sécantes
de $C\hookrightarrow\mathbb{P}(H^{o}(C,\Omega_{C})^{*})=\mathbb{P}(H^{o}(\Omega_{S})^{*})$
et la courbe $C\hookrightarrow\mathbb{P}^{g-1}$ est le lieu des points
exceptionnels. \\
Si $g\geq5$, on a: \[
\deg\psi=1,\,\deg F=2(g-1)(g-3),\]
où $\deg\psi$ est le degré de l'application cotangente et $\deg F$
celui de son image. Si $g=4$, alors $\deg\psi=6$.
\end{prop}
Rappelons que pour un point $s$ de $S$, nous avons noté $L_{s}$
l'image par l'application $\psi$ de $\pi^{-1}(s)$. Pour un point
$s$ de $S$, on note $\mathcal{I}_{s}$ l'idéal de définition de
$s$. Pour un diviseur $D$ de $C$, on note $\Omega_{C}(D)$ le fibré
$\Omega_{C}\otimes\mathcal{O}_{C}(D)$. 
\begin{lem}
\label{isomorphisme tilde}Supposons que $C$ soit non-hyperelliptique
de genre supérieur ou égal à $4$.\\
a) Le morphisme naturel $\vartheta:C^{(2)}\rightarrow J(C)$ est
un morphisme d'Albanese et un plongement. Le fibré cotangent $\Omega_{S}$
est engendré par ses sections globales.\\
b) Les invariants numériques de $S=C^{(2)}$ sont :\[
c_{1}^{2}[S]=(g-1)(4g-9),\, c_{2}[S]=(g-1)(2g-3),\, q=g.\]
c) Soit $s=P_{1}+P_{2}$ un point de $S$, l'image par l'isomorphisme
trace de l'espace $H^{o}(C,\Omega_{C}(-P_{1}-P_{2}))$ est $H^{o}(S,\mathcal{I}_{s}\Omega_{S})$.
\\
d) Soit $s=P_{1}+P_{2}\in C^{(2)}$ avec $P_{1}\not=P_{2}$. La
droite $L_{s}$ est la droite passant par les points $\phi(P_{1})$
et $\phi(P_{2})$.\end{lem}
\begin{proof}
Le point a) utilise le fait que $C$ n'est pas hyperelliptique et
l'appendice \cite{Mumford}. Le point b) est classique. Le point c)
résulte de \cite{Griffiths}. Par le corollaire \ref{diff=0000E9rentielle de l'apli. cot.},
le point d) est une autre formulation de c).
\end{proof}
Le degré de la variété des sécantes à une courbe lisse de genre $g$
et de degré $d$ dans un espace projectif de dimension $\geq4$ est
égal a \[
-g+(d-1)(d-2)/2\]
donc le degré de l'image de $\psi$ est $2(g-1)(g-3)$ si $g\geq5$.
Par la proposition \ref{degre de psi}, on a $\deg\psi\deg F=2(g-1)(g-3)$
donc $\deg\psi=1$.

Soit $P\in C$, pour tout élément $s$ de $P+C$, la droite $L_{s}$
(correspondant au point $\mathcal{G}(s)$) passe par $\phi(P)$. Le
point $\phi(P)$ est donc sommet d'un cône et $P+C$ est une courbe
non-ample. On vérifie que $(P+C)^{2}=1$.\\
La proposition \ref{proposition produit sym=0000E9trique} est
donc démontrée.

Les revêtements étales de surfaces isogènes à un produit donnent des
exemples différents de surfaces possédant une infinité de courbes
non-amples telles que $C^{2}>0$.

\curraddr{Graduate School of Mathematical Sciences, The University of Tokyo,
3-8-1 Komaba, Meguro, Tokyo, 153-8914 Japan.}

\email{roulleau@ms.u-tokyo.ac.jp}
\end{document}